\documentclass
{article}
\usepackage[T2A]{fontenc}
\usepackage[cp1251]{inputenc}
\usepackage[english,russian]{babel}
\usepackage[tbtags]{amsmath}
\usepackage{amsfonts,amssymb,mathrsfs,amscd}
\usepackage{wasysym}
\usepackage{verbatim}
\usepackage{eucal}
\usepackage{upgreek}
\usepackage{graphicx}
\usepackage{enumerate}
\let\No=\textnumero

\usepackage[hyper]{izv2e}
\JournalName{}
\numberwithin{equation}{section}

\theoremstyle{plain}
\newtheorem{theorem}{Теорема}
\newtheorem{maintheorem}{Основная теорема}
\newtheorem{lemma}{Лемма}
\newtheorem{propos}{Предложение}

\newtheorem{thF}{Теорема Фростмана}

\theoremstyle{definition}
\newtheorem{definition}{Определение}
\newtheorem{proof}{Доказательство}
\newtheorem{remark}{Замечание}
\newtheorem{example}{Пример}

\renewcommand{\leq}{\leqslant} 
\renewcommand{\geq}{\geqslant}
\newcommand{\RR}{\mathbb{R}} 
\newcommand{\CC}{\mathbb{C}} 
\newcommand{\NN}{\mathbb{N}}

\DeclareMathOperator{\supp}{{\sf supp}}
\DeclareMathOperator{\dd}{\,{\mathrm  d\!}}

\DeclareMathOperator{\mes}{mes}

\begin{document} 
\title{Мероморфные функции и разности субгармонических функций в интегралах и разностная характеристика Неванлинны. III. Оценки интегралов по фрактальным множествам через обхват и меру Хаусдорфа}


\author[B.\,N.~Khabibullin]{Б.\,Н.~Хабибуллин}
\address{Башкирский государственный университет}
\email{khabib-bulat@mail.ru}

\date{}
\udk{517.547.26 : 517.547.28 : 517.574 : 517.987.1}

 \maketitle

\begin{fulltext}

\begin{abstract} 
Пусть  $U\not\equiv \pm\infty$ --- $\delta$-субгар\-м\-о\-н\-и\-ч\-е\-ская  функция, в окрестности замкнутого  круга радиуса $R$
с центром в нуле. В предшествующих двух частях  нашей работы были получены общие и явные  оценки интеграла от положительной части радиальной максимальной характеристики роста  ${\mathsf M}_U(t):=\sup\bigl\{U(z)\bigm| |z|=r\bigr\}$  по возрастающей функции интегрирования $m$ на отрезке $[0,r]$ через разностную характеристику Неванлинны и величины, связанные с функцией интегрирования $m$. В третьей части работы  
оценки этих величин даются через $h$-обхват и  $h$-меру Хаусдорфа компакта $S\subset [0,r]$, для которого  функция интегрирования $m$ постоянна на каждой открытой компоненте связности дополнения $[0,r]\setminus S$. Отдельно рассмотрен случай $d$-мерной меры Хаусдорфа компакта  $S$.

Библиография: 10  названий 

Ключевые слова:  $\delta$-субгармоническая функция,  заряд  Рисса,  характеристика Неванлинны, обхват  Хаусдорфа, мера Хаусдорфа, $d$-мерная мера Хаусдорфа, фрактальное множество

\end{abstract}

\markright{Мероморфные функции и разности субгармонических функций  \dots. III}


\section{Результаты предшествующих двух частей работы}\label{s10}

{\it Интервал $I$\/} на расширенной	вещественной оси  $\overline \RR$ --- {\it связное\/} подмножество в $\overline \RR$ с {\it левым концом\/} $\inf I\in \overline \RR$ и с {\it правым концом\/} $\sup I\in \overline \RR$. 
{\it Отрезок\/} с {\it концами\/} $a\leq b\in \overline \RR$ --- это интервал $[a,b]:=\bigl\{x\in \overline \RR\bigm| a\leq x\leq b\bigr\}$. Кроме того,  $[a,b):=[a,b]\setminus b$ (соответственно  $(a,b]:=[a,b]\setminus a$) --- {\it открытый справа}
(соответственно {\it слева}) и {\it замкнутый слева\/} (соответственно {\it справа}) интервал, 
$(a,b):=[a,b)\cap (a,b]$ --- {\it открытый интервал.\/}  Через 
$D(r)$, $\overline  D(r)$ и $\partial \overline  D(r)$ 
обозначим соответственно {\it открытый\/} и  {\it замкнутый круги,\/} а также  {\it  окружность\/}  в $\CC$ {\it радиуса $r\in \overline \RR^+$ с центром в нуле.\/}   Для $R\in \overline \RR^+$ и функции $v\colon D(R)\to \overline \RR$ 
  \begin{equation}\label{{MC}M}
{\sf M}_v(r):=\sup\bigl\{v(re^{i\theta})\bigm| \theta \in [0,2\pi)\bigr\} , \quad 0\leq r<R,
\end{equation}
--- {\it максимальная характеристика роста функции $v$  на окружностях  $\partial \overline D(r)$,}
а ${\sf M}_v^+(r)$ --- её положительная часть. 
Каждой $\delta$-субгармонической функции   $U\not\equiv \pm \infty$ на замкнутом круге $\overline D(R)$ с положительной частью $U^+$ сопоставляется  {\it разностная характеристика Невалинны\/} 
\begin{equation}\label{rT}
{\boldsymbol T}_U(r,R):=
\frac{1}{2\pi}\int_0^{2\pi} U^+(re^{i\varphi})\dd \varphi 
+ \int_{r}^{R}\frac{\varDelta_U^-\bigl(\overline D(t)\bigr)}{t}\dd t,
\quad 0\leq r <R\in \RR^+,
\end{equation}
где $\varDelta_U^-$ --- {\it нижняя вариация заряда Рисса функции $U$ на $\overline D(R)$.} 

Сформулируем ключевое определение из \cite[определение 1]{Kha21I}.
\begin{definition}\label{DefhR}
{\it Возрастающей функции\/}   $m\colon [0,r]\to  \RR$ {\it полной вариации\/}
\begin{equation}\label{{hR}wm}
 {\tt M} :=m(r)-m(0) \in \RR^+
\end{equation}
с  {\it модулем непрерывности\/} $\omega_m\colon \RR^+\to \RR^+$, который  можно задать  как 
\begin{equation}\label{{hR}h}
\omega_m(t)\underset{t\in \RR^+}{=}\sup\bigl\{ m(x)-m(x')\bigm|x-x'\leq t, \, 0\leq x'\leq x\leq r  \bigr\}\overset{\eqref{{hR}wm}}{\subset} [0,{\tt M}],
\end{equation}
будем сопоставлять {\it диаметр стабилизации\/}
\begin{equation}\label{{hR}R}
{\sf d}_m:=\inf\bigl\{t\in \RR^+\bigm| {\omega}_{m}(t)= {\tt M}\bigr\}=\inf {\omega}_{m}^{-1}({\tt M})\leq r.
\end{equation}
Саму функция $m$ на $[0,r]$ с тем же обозначением $m$ часто будем рассматривать как {\it  возрастающую продолженную на $\RR$} постоянными  значениями $m(r)$ на луче $(r,+\infty)$ и постоянными значениями $m(0)$ на отрицательном луче $-\RR^+$, очевидно, {\it без увеличения полной вариации\/} ${\tt M}$. {\it Носителем  непостоянства\/} $\supp m'\subset [0,r]$  будем называть множество, состоящее из всех точек на $\RR$, в каждой окрестности которых продолженная функция  $m$ принимает хотя бы два различных значения. Так, если продолженная  функция $m$ дифференцируема, то, очевидно, её носитель непостоянства --- это носитель ей производной. 
\end{definition}

Все {\it интегралы\/} (Римана\,-- или  Лебега\,--\,){\it Стилтьеса\/}
с нижним пределом интегрирования    $a\in \RR$ и и верхним пределом интегрирования $b\geq a$  понимаем как интеграл по отрезку $[a,b]$, если не оговорено иное.

\begin{maintheorem}[{\rm \cite[основная теорема]{Kha21I}}]\label{th1} Пусть\/  $0< r<R\in \RR^+$. 
  
 Если  для модуля непрерывности $\omega_m$   возрастающей функции $m\colon [0,r]\to\RR$ выполнено условие Дини
\begin{equation}\label{{hR}i}
\int_0^{4R}\frac{{\omega}_{m}(t)}{t}\dd t<+\infty,
\end{equation}
то  для любой $\delta$-субгармонической функции $U\not\equiv\pm\infty$ на $\overline D(R)$ существует интеграл Лебега\,--\,Стилтьеса   с верхней оценкой 
\begin{equation}\label{U}
\int_0^r {\sf M}_{U}^+(t)\dd m(t)\leq \frac{6R}{R-r}
{\boldsymbol   T}_U(r,R) \max\biggl\{{\tt M}, \int_0^{{\sf d}_m}\ln \frac{4R}{t}\dd {\omega}_{m}(t)\biggr\}, 
\end{equation}
 где первый аргумент $r$ в ${\boldsymbol   T}_U(r,R)$ 
можно заменить на любое  $r_0\in [0,r]$, а   последний  интеграл Римана\,--\,Стилтьеса в \eqref{U} под операцией $\max$ --- на сумму
\begin{equation}\label{kint}
\int_0^{{\sf d}_m} \frac{{\omega}_{m}(t)}{t}\dd t+
{\tt M}\ln \frac{4R}{{\sf d}_m}\geq \int_0^{{\sf d}_m} \ln\frac{4R}{t}\dd {\omega}_{m}(t).
\end{equation} 
\end{maintheorem}

\begin{theorem}[{\rm \cite[теорема 1]{Kha21II}}]\label{th2}
Пусть $0<r\in \RR^+$, $h\colon [0,r]\to \RR^+$  --- непрерывная функция с  $h(0)=0$, дифференцируемая на $(0,r)$ и    удовлетворяющая условию   
\begin{equation}\label{{chrh}C}
{\sf s}_h:=\sup_{0< t< r} \frac{h(t)}{th'(t)}<+\infty.
\end{equation}
Тогда $h$ строго возрастает,  а  для любой возрастающей  функции $m\colon [0,r]\to\RR$ с полной вариацией  
${\tt M}\overset{\eqref{{hR}wm}}{:=}m(r)-m(0) \in \RR^+$ и 
модулем непрерывности 
\begin{equation}\label{{chrh}h}
{\omega}_m(t)\overset{\eqref{{hR}h}}{\leq} h(t)\text{ при всех $t\in [0,r]$}
\end{equation}
существует единственный прообраз  $h^{-1}({\tt M})\leq r$,  с которым   для  любой   $\delta$-суб\-г\-а\-р\-м\-о\-н\-и\-ч\-е\-с\-к\-ой  функции   $U\not\equiv \pm \infty$ на замкнутом круге  $\overline D(R)$ радиуса   $R>r$ существует  интеграл Лебега\,--\,Стилтьеса с верхней оценкой   
\begin{equation}\label{Uh}
\int_0^r {\sf M}_{U}^+(t)\dd m(t) \leq \frac{6R}{R-r}{\boldsymbol  T}_U( r, R)\,
 {\tt M}\ln \frac{4e^{{\sf s}_h}R}{h^{-1}({\tt M})}, 
\end{equation}
где первый аргумент $r$ в ${\boldsymbol   T}_U(r,R)$ 
можно заменить на любое число $r_0\in [0,r]$. 
\end{theorem}

\section{
Оценки интегралов через   $h$-меру и $h$-обхват Хаусдорфа} 

\subsection{Общий случай}
Для множества $X$, как обычно, $2^X$  --- {\it множество всех подмножеств  в $X$.} 
Следующие общие понятия  используются  в статье лишь  для интервалов на  $\RR$, поэтому далее 
  только такой случай и обсуждается.

\begin{definition}[{(\cite{Carleson}, \cite[\S~2.10]{Federer}, \cite{Rodgers}, \cite{HedbergAdams}, \cite{EG},
\cite{Eid07}, \cite{VolEid13})}]\label{defH}
Пусть $I\subset \RR$ --- интервал длины $r$, $h\colon [0, r]\to \overline \RR^+$ --- функция,  
$0<l\in \overline \RR^+$.  Функцию 
\begin{equation}\label{mr}
{\frak m}_h^{\text{\tiny $l$}}\colon S\underset{S\subset I}{\longmapsto}  \inf \Biggl\{\sum_{j\in N} h(b_j-a_j)
\biggm| N\subset \NN,\,  S\subset \bigcup_{j\in N} 
[a_j,b_j], \begin{cases}
[a_j,b_j]\subset I\\
0\leq b_j-a_j< l
\end{cases}\hspace{-3mm}\Biggr\} 
\end{equation}
на $2^I$ со значениями в $\overline \RR^+$ называем  {\it $h$-обхватом  Хаусдорфа\/} {\it диаметра $l$ на $I$}, не указывая  диаметр при $l=+\infty$, а  {\it $h$-обхватом  Хаусдорфа\/} {\it диаметра\/ $0$ на $I$\/} называем  функцию  
 \begin{equation}\label{hH0}
{\frak m}_h^{\text{\tiny $0$}}\colon S\underset{S\subset I}{\longmapsto}\sup_{l>0} {\frak m}_h^{\text{\tiny $l$}}(S) 
\end{equation} 
на $2^I$ со значениями в $\overline \RR^+$. При  $d\in \RR^+$ для   функции 
\begin{equation}\label{hd}
h_d\colon t\underset{t\in \RR^+}{\longmapsto} 
c_dt^d, \quad\text{где } 
c_d:=\dfrac{\pi^{d/2}}{2^d \Gamma(d/2+1)}, 
\quad  \Gamma\text{ \it  --- гамма-функция},
\end{equation}
$h_d$-обхват  Хаусдорфа  диаметра $l\in \overline \RR^+$ на $I$
называем  {\it $d$-мерным  обхватом Хаусдорфа диаметра $l$ на $I$\/}, который обозначаем как 
\begin{equation}\label{dmhl}
d\text{\tiny-}{\frak m}^{\text{\tiny $l$}}\overset{\eqref{hd}}{:=}{\frak m}_{h_d}^{\text{\tiny $l$}}, \quad d\in \RR^+, \quad l\in \overline \RR^+ .
\end{equation}
\end{definition}

Для множества $X$ функция $\mu \colon 2^X\to \overline \RR^+$ 
 {\it счётно субаддитивна на $X$,\/} если  {\it справедлива импликация\/}
\begin{equation}\label{subadd}
\biggl(N\subset \NN, \quad X_0\subset  \bigcup_{j\in N}X_j\subset X\biggr)
\Longrightarrow  \biggl( \mu(X_0) \leq \sum_{j\in N}\mu(X_j)\biggr).
\end{equation} 
Счётно субаддитивная функция  $\mu$ на $X$ 
  при    $\mu(\varnothing)=0$  называется  {\it внешней мерой  на  $X$.\/} Понятия {\it меры Бореля\/} и  {\it регулярной меры\/}  на  подмножествах в $\RR$ и $\CC$ общепринятые, а  $\supp \mu$ ---  {\it носитель\/} меры Бореля $\mu$.  {\it Мера Радона\/} --- регулярная мера Бореля, конечная на компактах.   

Тривиальные  свойства  понятий, введённых в  определении \ref{defH}, объединим в 
\begin{propos}\label{prhm} Для   $l\in \overline \RR^+$   и интервала   $I\subset \RR$ 
\begin{enumerate}[{\rm (i)}]
\item\label{hmi} $h$-обхваты Хаусдорфа ${\frak m}_h^{\text{\tiny $l$}}$ диаметра $l$ на $I$ счётно субаддитивны\/ \eqref{subadd};

\item\label{hmiI} значения    ${\frak m}_h^{\text{\tiny $l$}}(S)$ убывают по $l$ и  существует предел   
\begin{equation}\label{hH}
\lim_{0<l\to 0} {\frak m}_h^{\text{\tiny $l$}}(S)\overset{\eqref{hH0}}{=}{\frak m}_h^{\text{\tiny $0$}}(S)
\geq {\frak m}_h^{\text{\tiny $l$}}(S)\geq {\frak m}_h^{\text{\tiny $\infty$}}(S)
 \quad \text{\it для каждого  $S\subset I$}; 
\end{equation}
\item\label{hmiii} если $h(0)=0$, то все ${\frak m}_h^{\text{\tiny $l$}}$ ---  внешние  меры, а    ${\frak m}_h^{\text{\tiny $0$}}$  --- регулярная  мера Бореля на $I$, называемая  $h$-мерой  Хаусдорфа на $I$ и соответственно  
$d\text{\tiny-}{\frak m}^{\text{\tiny $0$}}$ --- $d$-мерная мера   Хаусдорфа;


\item\label{hmv} $0\text{\tiny-}{\frak m}^{\text{\tiny $0$}}(S)$ равно числу элементов в $S\subset I$;

\item\label{hmvi} если   $d>1$, то  $d$-мерная мера Хаусдорфа $d\text{\tiny-}{\frak m}^{\text{\tiny $r$}}$ нулевая на $I$.
\end{enumerate}
\end{propos}

Здесь и далее классические и широко известные свойства $h$-обхватов и $h$-мер Хаусдорфа из основных источников, указанных в определении  \ref{defH}, как и легко вытекающие из него, используются в адаптации именно для $\RR$ без явно прописанных  конкретных  ссылок или доказательств, что относится и к утверждениям предложения \ref{prhm}. Ввиду свойств   \eqref{hmv}--\eqref{hmvi} всюду  $d$-мерную меру Хаусдорфа на $\RR$ рассматриваем и используем только в нетривиальных случаях  $d\in (0,1]$. 

\begin{example} Использованная в \cite[Введение]{Kha21I} линейная мера Лебега  $\mes$ на $\RR$   --- это  в точности $1$-мерная  мера Хаусдорфа $1\text{\tiny-}{\frak m}^{\text{\tiny $0$}}$ на $\RR$, т.е. крайняя ситуация. 
\end{example} 

\begin{theorem}\label{th3} Если $m\colon [0,r]\to \RR$ --- возрастающая  функция
с носителем непостоянства\/ $\supp m'\subset S\subset [0,r]$,  $h\colon [0,r]\to \RR^+$ --- возрастающая функция с возрастающей  полунепрерывной сверху  регуляризацией $h^*$ на $[0,r]$,
и, подобно \eqref{{chrh}h}, выполнено неравенство 
 \begin{equation}\label{{chrh}h++}
h^*(t)\geq {\omega}_m(t) \quad\text{для каждого  $t\in [0,r]$},
\end{equation}   
то для  $h^*$-обхватов  Хаусдорфа 
на $I:=[0,r]$  выполнены   неравенства 
\begin{equation}\label{Mmm}
m(r)-m(0)\overset{\eqref{{hR}wm}}{=:}{\tt M}\leq {\frak m}_{h^*}^{\text{\tiny $l$}}(S) 
\leq h(r) 
\quad\text{при любых $l\in \overline \RR^+$},
\end{equation}
исходя из чего полную вариацию\/ ${\tt M}$  можно заменить 
\begin{enumerate}[{\rm (i)}]
\item\label{1Mm} в основной теореме в правой части \eqref{U} и  в левой части \eqref{kint}  на ${\omega}^*_m$-об\-х\-в\-ат Хаусдорфа ${\frak m}_{{\omega}^*_m}^{\text{\tiny $l$}}(S)$ диаметра $l\in \overline \RR^+$ на  $[0,r]$ 
множества $S$;
 
\item\label{2Mm} в неравенстве  \eqref{Uh} теоремы\/ {\rm \ref{th2}} на $h$-обхват Хаусдорфа ${\frak m}_{h}^{\text{\tiny $l$}}(S)$ диаметра $l\in \overline \RR^+$ на отрезке $[0,r]$ множества $S$ ввиду того, что 
\begin{equation}\label{inFp}
{\tt M}\ln \frac{4e^{{\sf s}_h}R}{h^{-1}({\tt M})}\leq {{\frak m}_{h}^{\text{\tiny $l$}}(S)}\ln \frac{4e^{{\sf s}_h}R}{h^{-1}({\frak m}_{h}^{\text{\tiny $l$}}(S))}.
\end{equation}
\end{enumerate}
\end{theorem}

Допускаемая теоремой \ref{th3}\eqref{2Mm} 
замена  в неравенстве  \eqref{Uh} теоремы\/ {\rm \ref{th2}}  выражения из левой части \eqref{inFp} на, вообще говоря, 
 большее выражение из правой части  \eqref{inFp} с  $h$-обхватом  Хаусдорфа ${\frak m}_{h}^{\text{\tiny $\infty$}}(S)$
может ослабить  оценку \eqref{Uh} в некотором смысле лишь на абсолютную постоянную-множитель, что отражает 

\begin{theorem}\label{th4}
Существует такая абсолютная постоянная $A\geq 1$, что для любого  $0<r\in \RR^+$, для всякого  
 компакта $S\subset [0,r]$ и  для каждой функции $h\colon [0,r]\to \RR^+$, удовлетворяющей всем условиям теоремы\/ {\rm \ref{th2}} с постоянной ${\sf s}_h\overset{\eqref{{chrh}C}}{>}0$,
  найдётся  возрастающая функция  $m$ на $[0,r]$ с полной вариацией ${\tt M}>0$,  с носителем непостоянства 
$\supp m' \subset S$, с модулем непрерывности, удовлетворяющим  \eqref{{chrh}h}, для которой одновременно  с неравенством \eqref{Uh} для произвольной $\delta$-субгармонической функции $U\not\equiv\pm\infty$ на шаре 
$\overline D(R)$ радиуса $R>r$ и неравенством \eqref{inFp}  выполнено противоположное к \eqref{inFp} неравенство с множителем $A$ перед\/ ${\tt M}$  вида
\begin{equation}\label{inFp-}
A\, {\tt M}\ln \frac{4e^{{\sf s}_h}R}{h^{-1}({\tt M})}\geq   {{\frak m}_{h}^{\text{\tiny $\infty$}}(S)}\ln \frac{4e^{{\sf s}_h}R}{h^{-1}({\frak m}_{h}^{\text{\tiny $\infty$}}(S))}.
\end{equation}
\end{theorem}

\subsection{Случай  $d$-мерного обхвата  Хаусдорфа диаметра $l\geq r$}

\begin{theorem}\label{cor1}
Пусть  $0<r\leq l\in \overline \RR^+$, $d\in (0,1]$. Для любой  возрастающей  функции $m\colon [0,r]\to\RR$ с носителем непостоянства 
 $\supp m'\subset S \subset [0,r]$ и модулем непрерывности, удовлетворяющим для некоторого $b\in \RR^+$ неравенству  
\begin{equation}\label{{chrh}hd}
{\omega}_m(t)\overset{\eqref{{hR}h}}{\leq} bt^d\quad \text{при  всех $t\in [0,r]$},
\end{equation}
и для  каждой   $\delta$-суб\-г\-а\-р\-м\-о\-н\-и\-ч\-е\-с\-к\-ой  функции   $U\not\equiv \pm \infty$ на  круге $\overline D(R)$
 радиуса $R>r$  существует  интеграл Лебега\,--\,Стилтьеса с верхней оценкой 
\begin{equation}\label{Uhd}
\int_0^{r} {\sf M}_{U}^+(t)\dd m(t) \leq \frac{24b}{d}
\frac{R}{R-r}
{\boldsymbol  T}_U(r,R)
\, d\text{\tiny-}{\frak m}^{\text{\tiny $l$}}(S)
\ln \frac{eR^d}{d\text{\tiny-}{\frak m}^{\text{\tiny $l$}}(S)}
\end{equation}  
где первый аргумент $r$ в ${\boldsymbol   T}_U(r,R)$ 
можно заменить на любое число $r_0\in [0,r]$. 
\end{theorem}
Эта теорема --- результат  пересечения теорем \ref{th2} и \ref{th3}, поэтому сразу можно дать 
\begin{proof}[теоремы \ref{cor1}] Положим  $h(t):=bt^d$ при $t\in [0,r]$, откуда  
\begin{equation}\label{h-1x}
{\sf s}_h\overset{\eqref{{chrh}C}}{=}\frac{1}{d}<+\infty, \quad h^{-1}(x)=\Bigl(\frac{x}{b}\Bigr)^{1/d}, \quad h\overset{\eqref{hd}}{=}\frac{b}{c_d}h_d,
\end{equation}
где  из определения \eqref{hd}  имеем  оценку снизу  для числа 
\begin{equation}\label{cd}
c_d\overset{\eqref{hd}}{:=}\dfrac{\pi^{d/2}}{2^d \Gamma(d/2+1)}\geq \frac{1}{2\Gamma(d/2+1)}\geq \frac{1}{4}
\quad\text{при $d\in (0,1]$.}
\end{equation}

По условию \eqref{{chrh}hd} выполнено и условие \eqref{{chrh}h} теоремы \ref{th2}, из которой 
\begin{equation}\label{Uhpr}
\int_0^r {\sf M}_{U}^+(t)\dd m(t) \leq \frac{6R}{R-r}
{\boldsymbol  T}_U( r, R)\, {\tt M}\ln \frac{4e^{1/d}R}{({\tt M}/b)^{1/d}}<+\infty.
\end{equation}
Кроме того, по неравенству \eqref{Mmm} теоремы \ref{th3} имеем 
\begin{equation}\label{Mdh}
{\tt M}\overset{\eqref{Mmm}}{\leq}{\frak m}_{h}^{\text{\tiny $l$}}(S)
\overset{\eqref{h-1x}}{=}\frac{b}{c_d}{\frak m}_{h_d}^{\text{\tiny $l$}}(S)
\overset{\eqref{dmhl}}{=}\frac{b}{c_d} d\text{\tiny-}{\frak m}^{\text{\tiny $l$}}(S)
\overset{\eqref{cd}}{\leq} 4b\,d\text{\tiny-}{\frak m}^{\text{\tiny $l$}}(S) \leq 4br^d,
\end{equation}
где последнее неравенство при  $r\leq l\in \overline \RR^+$
следует из определения \ref{defH} для $d$-мерного обхвата Хаусдорфа диаметра $l$ для множества $S$, содержащегося в отрезке $[0,r]$ длины не меньше $l$.

Для  выражения из \eqref{Uhpr}, содержащего ${\tt M}$, учитывая условие $0<d\leq 1$ имеем
\begin{equation}\label{MRd}
{\tt M}\ln \frac{4e^{1/d}R}{({\tt M}/b)^{1/d}}\leq {\tt M}\frac{1}{d}\ln \frac{e4^dR^db}{{\tt M}}
\leq \frac{1}{d}{\tt M}\ln \frac{e4bR^d}{{\tt M}}
\end{equation}
Для $0<B\in \RR^+$ функция $x\underset{x\in \RR^+}{\longmapsto} x\ln \frac{eB}{x}$ возрастающая на  $[0,B]$.
Если выбрать $B:=4bR^d$, то по  неравенствам  \eqref{Mdh}, где  $r<R$, в правой части \eqref{MRd}
можем заменить ${\tt M}$ на промежуточное $4b \,  d\text{\tiny-}{\frak m}^{\text{\tiny $l$}}(S)$ из  \eqref{Mdh}, 
что даёт
\begin{equation*}
{\tt M}\ln \frac{4e^{1/d}R}{({\tt M}/b)^{1/d}}
\leq \frac{1}{d}{4b \,  d\text{\tiny-}{\frak m}^{\text{\tiny $l$}}(S)}\ln \frac{e4bR^d}{{4b \,  d\text{\tiny-}{\frak m}^{\text{\tiny $l$}}(S)}}=\frac{4b}{d}{\,  d\text{\tiny-}{\frak m}^{\text{\tiny $l$}}(S)}\ln \frac{eR^d}{{d\text{\tiny-}{\frak m}^{\text{\tiny $l$}}(S)}}.
\end{equation*}
Последнее  согласно \eqref{Uhpr} влечёт за собой  неравенство \eqref{Uhd} теоремы \ref{cor1}.
\end{proof}
\begin{remark} В случае $d=1$ и $1\text{\tiny-}{\frak m}^{\text{\tiny $0$}}$ меры Хаусдорфа, совпадающей с линейной мерой  Лебега $\mes$, в качестве возрастающей функции $m$, по которой  интегрируются максимальные характеристики роста ${\sf M}_U^+$, всегда можно выбирать  функции распределения  сужений меры $\mes$ на произвольные  борелевские подмножества на отрезках из $\RR^+$, поскольку все такие функции распределений имеют модуль непрерывности,  не превышающий тождественной функции, и поэтому  удовлетворяют условию Дини. 
Но при $d\in (0,1)$, вообще говоря,   функции распределения  сужений $d$-мерной меры Хаусдорфа 
 на  борелевские подмножества конечной $d\text{\tiny-}{\frak m}^{\text{\tiny $0$}}$-меры на отрезках из $\RR^+$ совсем не обязаны обладать модулем непрерывности, удовлетворяющим условию Дини. 
\end{remark}

\section{Один фундаментальный результат теории потенциала}
В настоящей статье только на $\RR$ формулируется и используется    

\begin{thF}[\hspace{-1mm}{\rm (\cite[II, теорема 1]{Carleson}, \cite[теорема 5.1.12]{HedbergAdams}, \cite[теорема 2.11]{AE})}]  Справедливы следующие два утверждения:
\begin{enumerate}[{\rm I.}]
\item\label{IF} Пусть $I\subset \RR$ --- интервал длины $r$,   функция $\mu\colon 2^{I}\to \overline \RR^+$ счётно субаддитивна и для  функции  $h\colon [0,r]\to \overline \RR^+$ имеют место  неравенства  
\begin{equation}\label{mur}
\mu\bigl([a, b]\bigr)\leq h(b-a)\quad\text{для всех отрезков $[a,b]\subset I$.}
\end{equation}
Тогда для $h$-обхватов Хаусдорфа любого диаметра $l\in \overline \RR^+$ на $I$ имеем
\begin{equation}\label{muS}
{\frak m}_{h}^{\text{\tiny $l$}}(S)\geq \mu(S)  \quad\text{при  всех $S\subset I$.} 
\end{equation}

\item\label{IIF} Существует такое  положительное число  $A$, что для каждой возрастающей функции $h\colon \RR^+\to \RR^+$  с $h(0)=0$ и  для любого компакта $E$ в $\RR$ найдётся  мера Радона 
$\mu\neq 0$ на $\RR$ с носителем $\supp \mu \subset E$ со свойствами \eqref{mur} для $I:=\RR$, а также \eqref{muS} для всех борелевских подмножеств $S\subset \RR$, но при этом ещё и 
\begin{equation}\label{muhr}
 {\frak m}_h^{\text{\tiny $\infty$}}(E) \leq A\mu(E).
\end{equation} 
\end{enumerate}
\end{thF}
В известных нам  формулировках теоремы Фростмана обе части \ref{IF} и \ref{IIF}
 даются 
с едиными посылками и, как следствие, с перегрузкой условий на функции $h$ и $\mu$ в тривиальной первой части \ref{IF}.  Поэтому приведём 

\begin{proof}[части \ref{IF}]
Пусть, как в  \eqref{mr} из определения \ref{defH} с  $l\in \overline \RR^+\setminus 0$, 
\begin{equation}\label{SBj}
S\subset \bigcup_{j\in N} [a_j,b_j], \quad N\subset \NN, \quad  [a_j,b_j]\subset I, \quad  0\leq b_j-a_j<l .
\end{equation}
Тогда из счётной аддитивности функции $\mu$ получаем 
\begin{equation*}
\mu(S)\overset{\eqref{subadd}}{\leq} \sum_{j\in N}\mu\bigl( [a_j,b_j]\bigr)\overset{\eqref{mur}}{\leq} 
\sum_{j\in N}h( b_j-a_j),
\end{equation*}
и применение к крайним частям неравенств  точной нижней грани в ограничениях   \eqref{SBj}, от которых левая часть не зависит, сразу устанавливает \eqref{muS} для любого $l\in \overline \RR^+\setminus 0$, откуда по определению  \eqref{hH0} получаем \eqref{muS} и для  $l=0$. Часть \ref{IF} доказана.
\end{proof}

\section{Доказательство теоремы \ref{th3}}
 Выведем сначала части \eqref{1Mm}--\eqref{2Mm} из неравенств \eqref{Mmm}.
Часть \eqref{1Mm} сразу следует из   \eqref{Mmm} при $h:={\omega}_m\leq {\omega}_m^*$, так как полная вариация  ${\tt M}$ входит в \eqref{U} и \eqref{kint} в форме, позволяющей заменить  ${\tt M}$ на любое большее число. 

Для вывода части \eqref{2Mm} потребуется элементарная 
\begin{lemma}\label{lem2} Пусть  функция $h$ такая же, как в теореме\/ {\rm \ref{th2}},  с числом  ${\sf s}_h>0$ 
из \eqref{{chrh}C}, а значит строго возрастающая на $[0,r]$, и  пусть $r\leq B\in \RR^+$.  Тогда
\begin{equation}\label{h-1}
x\longmapsto x\ln \frac{Be^{{\sf s}_h}}{h^{-1}(x)}, \quad x\in \bigl[0, h(r)\bigr] ,
\end{equation}
--- возрастающая функция на отрезке  $\bigl[0, h(r)\bigr]$.
\end{lemma}
\begin{proof}[леммы \ref{lem2}] Произведём замену $t:=h^{-1}(x)\in [0,r]$ и перейдём от функции 
\eqref{h-1} к функции 
\begin{equation}\label{h-1-}
t\longmapsto h(t)\ln \frac{Be^{{\sf s}_h}}{t}, \quad t\in [0, r]. 
\end{equation}
Ввиду строгого возрастания непрерывной функции $h$  на $[0,r]$ достаточно показать, что возрастает функция  \eqref{h-1-}. Дифференцирование этой функции на открытом интервале $(0,r)$ даёт 
\begin{equation*}
h'(t)\ln \frac{Be^{{\sf s}_h}}{t}-\frac{h(t)}{t}\overset{\eqref{{chrh}C}}{\geq} 
h'(t)\ln \frac{Be^{{\sf s}_h}}{t}-{\sf s}_hh'(t)=h'(t)\ln \frac{B}{t}\geq 0
\end{equation*} 
на $(0,r)$ при $B\geq r$, откуда следует возрастание функции $h$ на открытом интервале  $(0,r)$, а в силу непрерывности и на отрезке $[0,r]$. 

Лемма \ref{lem2} доказана. 
\end{proof}

В неравенстве  \eqref{Uh} теоремы \ref{th2} полная вариация ${\tt M}$
встречается дважды исключительно  в составе  выражения
$$
{\tt M}\ln \frac{Be^{{\sf s}_h}}{h^{-1}({\tt M})}, \quad\text{где  $B:=4R\geq r$} 
$$
имеющего вид \eqref{h-1}, и  по лемме \ref{lem2} и неравенству \eqref{Mmm} замена ${\tt M}$ на 
${\frak m}_{h^*}^{\text{\tiny $l$}}(S)\leq h(r)$ не уменьшит это выражение. При этом в силу непрерывности $h$
в теореме \ref{th2}
имеем $h^*=h$ и  ${\frak m}_{h^*}^{\text{\tiny $l$}}={\frak m}_{h}^{\text{\tiny $l$}}$.

По определениям \eqref{mr}--\eqref{hH0} ввиду $S\subset [0,r]$, очевидно,  выполнены соотношения   ${\frak m}_{h^*}^{\text{\tiny $l$}}(S)\leq h^*(r-0)=h(r)$ для всех $l\in \overline \RR^+$ и  $S\subset [0,r]$, что доказывает второе  неравенство в  \eqref{Mmm}.

Для доказательства оставшегося первого неравенства в \eqref{Mmm} напомним, что {\it меру Лебега\,--\,Стилтьеса\/} $\mu$ на $[0,r]$ для {\it возрастающей\/} функции $m$ на отрезке $[0,r]\subset \RR^+$ можно сначала определить  
 через возрастающее продолжение её на $\RR$ без увеличения полной вариации ${\tt M}$, как в определении \ref{DefhR},  на   открытых слева и замкнутых справа интервалах $(a,b]\subset \RR$  равенствами 
\begin{equation}\label{m*}
 \mu\bigl((a,b]\bigr):=\lim\limits_{b<x\to b} m(x)-\lim\limits_{a<x\to a} m(x),
\end{equation}
 с последующим лебеговским продолжением на подмножества в $\RR$. По построению {\it носитель\/} $\supp \mu$ меры Лебега\,--\,Стилтьеса $\mu$ лежит на  $[0,r]$ и  совпадает с носителем непостоянства $\supp m'$. 
Из \eqref{m*} по определению \eqref{{hR}h} модуля непрерывности  ${\omega}_m$ получаем 
\begin{equation*}
 \mu\bigl([a,b]\bigr)\overset{\eqref{{hR}h}}{\leq}\lim_{b-a<t\to b-a} {\omega}_m(t)
\end{equation*}
откуда  по условию \eqref{{chrh}h++} для  полунепрерывной сверху регуляризации
$$
\mu\bigl([a,b]\bigr) \overset{\eqref{{chrh}h++}}{\leq} \lim_{b-a\leq t\to b-a}  h^*(t)=h^*(b-a)
\quad\text{при $[a,b]\subset [0,r]$}
$$
ввиду полунепрерывности сверху и возрастания функции $h^*$ на $[0,r]$. Последнее означает, что выполнено 
условие \eqref{mur} из части \ref{IF} теоремы Фростмана с полунепрерывной сверху регуляризацией $h^*$ в роли $h$, и  
по заключению \eqref{muS} имеем неравенство  
$\mu(S)\leq {\frak m}_{h^*}^{\text{\tiny $l$}}(S)$ при каждом  $ S\subset [0,r]$. Если $S\supset \supp \mu$, то левая часть этого неравенства равна полной $\mu$-мере отрезка $[0,r]\supset \supp \mu$, которая по построению 
\eqref{m*} равна   разности   $m(r)-m(0)\overset{\eqref{{hR}wm}}{=}{\tt M}$.

Теорема \ref{th3} доказана. 

\section{Доказательство теоремы \ref{th4}}
По заданным компакту $E:=S\subset [0,r]$ и функции $h$ из условия теоремы \ref{th2} с возрастающим продолжением на  $\RR$ с $h(t)\equiv h(r)$ при всех  $t\geq r$, как в определении \ref{DefhR},  по части \ref{IIF} теоремы Фростмана выберем меру Радона $\mu\neq 0$ с всеми прописанными в этой части  \ref{IIF} свойствами. В качестве возрастающей функции $m\colon [0,r]\to \RR^+$ выберем функцию распределения на $[0,r]$ это меры $m(t):=\mu \bigl([0,t]\bigr)$ при $t\in [0,r]$. 
По построению $\supp m'=\supp \mu \subset S\subset [0,r]$ и $\mu(S)=m(r)-m(0)=:{\tt M}$, свойство   \eqref{mur} по определению модуля непрерывности \eqref{{hR}h} означает, что $\omega_{m}(t)\leq h(t)$ при $t\in [0,r]$. Следовательно,  выполнено условие  \eqref{{chrh}h}  и по теореме \ref{th2} для  любой   $\delta$-суб\-г\-а\-р\-м\-о\-н\-и\-ч\-е\-с\-к\-ой  функции   $U\not\equiv \pm \infty$ на  $\overline D(R)$ радиуса   $R>r$ имеем \eqref{Uh}. При этом  из свойств \eqref{muS} и  \eqref{muhr} следует 
$$ 
{\tt M}=\mu(S)\leq {\frak m}_h^{\text{\tiny $\infty$}}(S) \leq A\mu(S)=A\,{\tt M}
$$
откуда 
$$
 {{\frak m}_{h}^{\text{\tiny $\infty$}}(S)}\ln \frac{4e^{{\sf s}_h}R}{h^{-1}({\frak m}_{h}^{\text{\tiny $\infty$}}(S))}
\leq A\,{\tt M}\ln \frac{4e^{{\sf s}_h}R}{h^{-1}({\frak m}_{h}^{\text{\tiny $\infty$}}(S))}
\leq A\,{\tt M}\ln \frac{4e^{{\sf s}_h}R}{h^{-1}({\tt M})},
$$
ввиду возрастания $h^{-1}$ в знаменателе. Это доказывает \eqref{inFp-} и теорему \ref{th4}.

\end{fulltext}

\end{document}